\documentclass[leqno,a4paper, 11pt]{article}
\usepackage{times}
\usepackage{amsmath}
\usepackage{amssymb}
\usepackage{amsfonts}

\newcommand{\hp}{\hspace*{\parindent}}

\usepackage{latexsym}
\newcommand{\hem}{\hspace*{1em}}

\newcommand{\hfl}{\hspace*{\fill}}
\newcommand{\nhp}{\hspace*{-\parindent}}

\newlength{\dede}     

\newcommand{\I}{\'{\i}}

\newcommand{\cao}{\c c\~ao}
\newcommand{\coes}{\c c\~oes}

\def\cuad{{\hfl \raggedleft\rule{1ex}{1ex}}}

\newtheorem{Teo}{Theorem}[section]
\newtheorem{Df}[Teo]{Definition}
\newtheorem{remark}[Teo]{Remark}
\newtheorem{Pro}[Teo]{Proposition}

\newtheorem{Lem}[Teo]{Lemma}
\newtheorem{notation}[Teo]{Notation}
\newtheorem{Cor}[Teo]{Corollary}

\newcommand{\n}{\noindent}
\newcommand{\dem}{\n{\bf Proof}}

\newcommand{\bc}{\begin{center}}
\newcommand{\ec}{\end{center}}

\newcommand{\vtres}{\vspace*{0.3cm}}

\newcommand{\Le}{\mathcal{L}}

\begin{document}
\date{}
\title{Distributive abstract logics and the Esakia duality}

\author{Andreas B.M.~Brunner\\Departamento de Matem\'atica\\Instituto de Matem\'atica\\Universidade Federal da Bahia - UFBA\\40170-110 Salvador - BA\\Brazil\\e-mail: andreas@dcc.ufba.br\and Darllan Concei\cao\ Pinto\\Departamento de Matem‡tica\\Instituto de Matem\'atica e Estat\I stica\\Universidade de S\~ao Paulo - USP\\- S\~ao Paulo - SP\\Brazil\\e-mail: darllan@ime.usp.com }

\maketitle

\begin{abstract}
In this paper we develop an almost general process to switch from abstract logics in the sense of Brown and Suszko to lattices. With this method we can establish dualities between some categories of abstract logics to the correspondent topological space categories. In more detail we will explain the duality between the category of abstract intuitionistic logics with intuitionistic morphisms and the category of Esakia spaces with the Esakia morphisms.
\end{abstract}

\section{Introduction}

\hp Abstract logics were introduced by earlier works from Brown and Suszko, cf \cite{brosus}. Basically, these authors see an abstract logic, as an intersection structure with or without greatest element. It is well known that this notion is in bijective correspondence with complete lattices, and also with closure operators in Tarski's sense.  It is also known, that algebraic intersection structures, algebraic lattices and compact Tarski-operators are in bijection. That is, every complete (algebraic) lattice gives us an (compact) Tarski-operator, and vice-versa. For details of this affirmation, we refer the reader  to \cite{bursan} or \cite{davpri}.  So defining an abstract logic as an intersection structure, we are able to ask some questions about them. For example we can work with intersection structures that have some more properties, which are introduced by the existence of some  connectives, see \cite{lewbru} and also \cite{fonver}. In the article of Bloom and Brown, cf \cite{blobro}, the authors work with abstract classical logics in a Boolean sense, but also abstract logics in a non-classical sense can be defined and worked with, cf. \cite{BL11}.

From that time on, many researches were made in this topic, between them also by the first author in joint work with S. Lewitzka, see \cite{BL11,  lew2, lewbru}.  The principal idea of this work is on one hand to establish an almost general method to {\sl switch} from abstract logics to lattices, and so to be able to generalize some duality results of the corresponding categories. Even almost in an easy manner we can go from abstract logics to lattices,  it is not immediately clear that the duality results, will hold. This is so, because the categories always carry with them morphisms, and in the beginning it is not clear that the kind of distinct logic maps {\sl do} give in fact on the other side the desired morphisms in the category considered. For example, it is known that in Stone's duality for distributive lattices, the category of distributive lattices with lattice morphisms is dually equivalent to the category of the spectral spaces, with the spectral functions as morphisms, cf. \cite{mir}. Also, for this category of distributive lattices there exist some dualities of bitopological nature, see for example  \cite{bezgabkur} and there is the well known Priestley duality, cf \cite{pri}. We will show that these results will hold also for our abstract distributive logics with the stable logic maps as introduced in \cite{lew2}. Clearly, a generalization is easily obtained for the Boolean abstract logics, Boolean algebras and Boolean spaces. In \cite{BL11} the authors establish a duality for the categories of intuitionistic and distributive abstract logics, with stable logic maps, and the categories of spectral spaces, with and without implication. These results are obtained using another strategy, and we think that also these results can be obtained by the method introduced here.

The paper is structured in the following way. We resume some important preliminaries about abstract logics, in the first section. In the second section, we explain the rather simple method of switching from abstract logics to lattices, and resume some duality results. Then in the last section we  will show in detail that the category of Heyting algebras with the Heyting morphism is indeed dually equivalent to the category of intuitionistic abstract logics with  intuitionistic logic maps. By this result, we obtain immediately that Esakia duality is valid for our intuitionistic abstract logics.

\section{The concepts of abstract logics}

\hp In this section we recall some definitions and results from abstract logics which are essentially given in the articles \cite{lewbru} and \cite{BL11}. For a more detailed presentation we refer the reader to these papers.

\begin{Df}\label{20}
An abstract logic $\mathcal{L}$ is given by $\mathcal{L}=(Expr_\mathcal{L},Th_\mathcal{L},\mathcal{C}_\mathcal{L})$, where $Expr_\mathcal{L}$ is a set of expressions (or formulas) and $Th_\mathcal{L}$ is a non-empty subset of the power set of $Expr_\mathcal{L}$, called the set of theories, such that the following intersection axiom is satisfied:

\nhp \hfl
$\text{If }\mathcal{T}\subseteq Th_\mathcal{L}\text{ and }\mathcal{T}\neq\emptyset,\text{ then }\bigcap\mathcal{T}\in
Th_\mathcal{L}. $\hfl

\nhp Furthermore, $\mathcal{C}_\mathcal{L}$ is a set of operations on $Expr_\mathcal{L}$, called (abstract) connectives.

\nhp (a) An abstract logic $\mathcal{L}$ is called {\bf regular} iff $Expr_\mathcal{L}$ is not a theory, i.e., $Expr_\mathcal{L}\notin Th_\mathcal{L}$. Otherwise, $\mathcal{L}$ is {\bf singular}. \\
(b)  A subset $A\subseteq Expr_\mathcal{L}$ is called {\bf consistent} iff $A$ is contained in some theory $T\in Th_\mathcal{L}$.\\
(c) A theory $T\in Th_\mathcal{L}$ is called {\bf $\kappa$-prime} ($\kappa\ge\omega$ a cardinal) iff for every non-empty set $\mathcal{T}\subseteq Th_\mathcal{L}$ of size $<\kappa$, $T=\bigcap\mathcal{T}$ implies $T\in\mathcal{T}$. In the case, in which $T$ is $\omega$-prime, we say that $T$ is {\bf prime}. A  {\bf totally} (or {\bf completely}) {\bf prime theory} is a theory which is $\kappa$-prime for all cardinals $\kappa\le\omega$. \\
(d) A set of theories $\mathcal{G}\subseteq Th_\mathcal{L}$ is called a {\bf generator set} for the logic $\mathcal{L}$ iff each theory is the intersection of some non-empty subset of $\mathcal{G}$. In the case, a minimal generator set exists, we say that $\mathcal{L}$ is {\bf minimally generated}.\\
(e) A theory $M \in Th_\mathcal{L}$ is called {\bf maximal} in a regular logic iff for every theory $T \in Th_\Le$ such that $M \subseteq T$, we have that $M = T$. \\
(f) An abstract logic  $\mathcal{L}$ is {\bf closed under union of chains} iff for any ordinal $\alpha>0$ and any chain of theories $\{T_i\mid i<\alpha\}$, the set $\bigcup_{i<\alpha}T_i$ is a theory. \\
(g) An abstract logic $\Le$ has a {\bf $\kappa$-disjunction}, $\bigvee$, iff for all sets of expressions $A \subseteq Expr_\Le$ of cardinality $ < \kappa$, all $T$  totally prime we have that:\\
\nhp \hfl $A \cap T \not = \emptyset$ \hem iff \hem $\bigvee A \in T$. \hfl

\end{Df}

Clearly, abstract logics have Tarski-consequence operators satisfying the three Tarski axioms. We can introduce them in the known way.

\begin{Df} \label{15}
Let $\Le$ be an abstract logic as in definition \ref{20} and $A \cup \{ a \} \subseteq Expr_\Le$.

\nhp (a) The consequence relation $\Vdash_\mathcal{L}$ is defined in the following way:

\nhp \hfl$A\Vdash_\mathcal{L} a$ \hem iff \hem  $ a\in\bigcap\{T\in Th_\mathcal{L}\mid A\subseteq T\}$. \hfl

\nhp (b) The consequence relation is called {\bf compact} or equivalently {\bf finitary} \hem iff \hem $A\Vdash_\mathcal{L} a$, implies the existence of a finite $A'\subseteq A$ such that $A'\Vdash_\mathcal{L}a$. \\
(c) The abstract logic $\mathcal{L}$ is called {\bf compact} iff every inconsistent set of formulas has a finite inconsistent subset.  \\
(d) The formula $a$ is {\bf valid} iff $a \in T$, for all theories $T \in Th_\Le$.
\end{Df}

Note that the notion of generator set corresponds to the concept of meet-dense subset of a meet-semilattice.

\begin{Pro}[cf. \cite{lewbru}]\label{30}
Let $\mathcal{L}$ be an abstract logic. Then we have the following:

\nhp (a) A set of expressions $T\subseteq Expr_\mathcal{L}$ is a theory iff $T$ is consistent and closed under $\Vdash_\mathcal{L}$ (i.e. $T$ is contained in some theory, and $T\Vdash_\mathcal{L}a$ implies $a\in T$). \\
(b) $\mathcal{L}$ is closed under union of chains (and regular) iff the consequence relation is compact (and there is a finite inconsistent set of formulas). \cuad

\end{Pro}

The first statement of \ref{30} follows easily from the definitions. The second statement follows from 2.17 in \cite{lewbru}, if $\mathcal{L}$ is regular. In the singular case, it follows from basic results about closure spaces (see for example, \cite{davpri}).

In \cite{lewbru}, it was proved in theorem (2.11.) that an abstract logic closed by union of chains is in fact minimally generated. This minimally set of generators was shown to be the totally prime theories $TPTh_{\Le}$. The proof there, was based on the well-ordering theorem and used methods of set theory. We want to give in the following a new proof of this theorem, using Zorn's Lemma and algebra.

 \begin{Teo} \label{MinGer}

Let $\mathcal{L}$ be an abstract regular logic closed by union of chains.  Then $\mathcal{L}$ is minimally generated by the set  $TPTh_{\Le}$.
\end{Teo}

\dem :  By hypothesis, $\mathcal{L}$ is closed by union of chains and so we have for every chain of theories $\mathcal{C} \subseteq Th_\Le$,  $\bigcup \mathcal{C} \in Th_\Le$. Let $T_0 \in Th_\Le$ arbitrary. We will show that $T_0$ is generated by totally prime theories.

The fact that $\mathcal{L}$ is a regular logic implies that  $T_0 \not = Expr_{\Le}$. For this reason we have $a \in Expr_{\Le}$ such that $a \not \in T_0$. Consider the following set,\\
\nhp \hfl $\mathcal{F} := \{ T \in Th_{\Le} | \hem T_0 \subseteq T \hem \& \hem a \not \in T \}$. \hfl \\
 It is clear that  $\mathcal{F}$ is not empty. Also $\mathcal{F}$ is partially ordered by inclusion, $\subseteq$. By hypothesis, for every chain $\mathcal{C} \subseteq \mathcal{F}$, $\bigcup \mathcal{C}$ is a theory. Because  $a \not \in \bigcup \mathcal{C}$, $\bigcup \mathcal{C}$ is an upper bound of  $\mathcal{C}$ and  Zorn«s Lemma can be applied. Denote by  $T_a$ a maximal element in $\mathcal{F}$. We show the following

{\bf Fact:} $T_a$ is totally prime. \\
Proof of fact:  Suppose that this is not so, i.e., $T_a$ is not totally prime. Then there exists a cardinal $\kappa \geq \omega$ and a family of theories of cardinality $\kappa$, say $\mathcal{T}_\kappa$ such that $T_a = \bigcap \mathcal{T}_\kappa$ and  $T_a$ is different of any element of  $\mathcal{T}_\kappa$, i.e., $T_a \subsetneq T$, $\forall T \in \mathcal{T}_\kappa$. \\
From the fact that  $T_a$ is  maximal with the property of being a theory which does not contain the formula  $a$, we must have for every  $T \in \mathcal{T}_\kappa$, $a \in T$. Observe now that
$(T_a \cup \{ a \})^{\Vdash_{\mathcal{L}}}$ is the least theory containing $T_a$ such that  $a \in (T_a \cup \{ a \})^{\Vdash_{\mathcal{L}}}$. For this, we have that $(T_a \cup \{ a \})^{\Vdash_{\mathcal{L}}} \subseteq T$, for every theory $T \in \mathcal{T}_\kappa$. Thus, $T_a$ is not an intersection of proper theories, and consequently, $T_a$ has to be totally prime.

Repeating this argument for all elements $b \not \in T_0$, we always obtain a totally prime theory $T_b$. So, $T_0 = \bigcap_{b \not \in T_0} T_b$. To see this equality, remark that always, if $a \in T_0$, then $a \in T_b$, for all $b \not \in T_0$. On the other hand, if $a \not \in T_0$, by construction $a \not \in T_a$ and thus, $a \not \in \bigcap_{b \not \in T_0} T_b$, finishing proof of theorem. \cuad

\begin{notation}
Let $MTh_\mathcal{L},TPTh_\mathcal{L},PTh_\mathcal{L}$ denote the sets of maximal, totally prime and prime theories of the abstract logic $\mathcal{L}$, respectively. It follows that $MTh_\mathcal{L}\subseteq TPTh_\mathcal{L}\subseteq PTh_\mathcal{L}$. Furthermore, $TPTh_\mathcal{L}$ is contained in any generator set. Thus, in a minimally generated logic $\mathcal{L}$, $TPTh_\mathcal{L}$ is the minimal generator set.
\end{notation}

The definition of \textit{intuitionistic abstract logic}, where the connectives are characterized by means of conditions over the minimal generator set, is given in \cite{lew2,lewbru}. We consider here also the notion of \textit{(bounded) distributive} abstract logic, and repeat the important definitions.

\begin{Df}\label{def1}
Let $\mathcal{L}=(Expr_\mathcal{L},Th_\mathcal{L},\mathcal{C}_\mathcal{L})$ be an abstract logic closed under union of chains. For a set $\{\vee,\wedge,\sim,\rightarrow\}$ of operators consider the following conditions. For all $a,b\in Expr_\mathcal{L}$ and for all $T\in TPTh_\mathcal{L}$:

\nhp (a) \hp $a\vee b\in T$ \hem iff \hem $ a\in T$ or $b\in T$.\\
(b) \hp $a\wedge b\in T$ \hem iff \hem $ a\in T$ and $b\in T$.\\
(c) \hp $\sim a\in T$ \hem iff \hem $ T\cup\{a\}$ is inconsistent. \\
(d) \hp $a\rightarrow b\in T$ \hem iff \hem  for all totally prime $T'\supseteq T$, if $a\in T'$ then $b\in T'$. \\
(e) \hp There is a formula $\top\in Expr_\mathcal{L}$ which is contained in every (totally prime) theory (i.e. $\top$ is valid). \\
(f) \hp There is a formula $\bot\in Expr_\mathcal{L}$ which is contained in no (totally prime) theory (i.e. $\bot$ is inconsistent).

\nhp If $\{\vee,\wedge\}\subseteq \mathcal{C}_\mathcal{L}$ and (a),(b) hold, then $\mathcal{L}$ is called a {\bf distributive abstract logic}. $\mathcal{L}$ is said to be {\bf bounded} iff in addition (e) and (f) hold. If $\mathcal{C}_\mathcal{L}=\{\vee,\wedge,\sim,\rightarrow\}$ and (a)-(d) hold, then $\mathcal{L}$ is an {\bf intuitionistic abstract logic}. An intuitionistic abstract logic $\mathcal{L}$ with $MTh_\mathcal{L}=TPTh_\mathcal{L}$ is called a {\bf classical} (or {\bf a boolean}) {\bf abstract logic}.
\end{Df}

We can show that an intuitionistic abstract logic is indeed  bounded.


In intuitionistic abstract logics the sets of maximal, totally prime and prime theories are in general distinct (see the discussion in \cite{lewbru}); these sets coincide in the classical case. In \cite{lewbru} we asked for a greatest set $\mathcal{T}\subseteq Th_\mathcal{L}$ of theories such that the conditions (a)-(d) of Definition \ref{def1} remain true if we replace $TPTh_\mathcal{L}$ by $\mathcal{T}$. We call such a set \textit{the set of complete theories} $CTh_\mathcal{L}$.  It was proved in \cite{lewbru} that $CTh_\mathcal{L}$ exists --- it is exactly the set of prime theories: $CTh_\mathcal{L}=PTh_\mathcal{L}$. In effect, it was shown a more general result considering appropriate notions of $\kappa$-disjunction and $\kappa$-conjuntion. Theorem 3.4 in \cite{lewbru} shows that in the presence of $\kappa$-disjunction, $CTh_\mathcal{L}$ is the set of all $\kappa$-prime theories --- this holds independently from the presence or absence of the other intuitionistic connectives.
In the case $\kappa=\omega$, this shows in particular that our notion of prime theory, introduced in an order-theoretic way, coincides with the usual notion of a prime theory $T$ in intuitionistic logic: $a\vee b\in T$ iff $a\in T$ or $b\in T$, for any formulas $a,b$.




For future use we will prove the following result. The proof is rather simple, but nonetheless we will elaborate it.

\begin{Lem} \label{Lemprelim}
Let $\Le$ an abstract intuitionistic logic and let $a, b \in Expr_\Le$ given. For every $T \in Th_\Le$ we have:\\
(a) \hp $a,b \in T$ \hem iff \hem $(a \wedge b) \in T$.\\
(b) \hp If $a \in T$ or $b \in T$ \hem then \hem $(a \vee b) \in T$. The implication the other way round is valid only for prime theories.\\
(c) \hp If $a \in T$ and $(a \to b) \in T$ \hem then \hem $b \in T$. (the theories are closed by modus ponens). \\
(d) \hp Let $P \in PTh_\Le$ a prime theory, then\\
\nhp \hfl $ a\rightarrow b\in P \hem \Leftrightarrow$ \hem  for every prime theory $Q\supseteq P,\ se\ a\in Q$ then $b\in Q$.   \hfl

\end{Lem}

\dem : Let  $T \in Th_\Le$ de a theory. As $\Le$ is minimally generated by $TPTh_\Le$, there exists $\tau \subseteq TPTh_\Le$ such that $T = \bigcap \tau$. Now we will show (a). Let $a, b \in Expr_\Le$ such that $a, b \in T$, i.e.,  $a, b \in Q$, for every $Q \in \tau$. By definition \ref{def1}, this is equivalent with  $(a \wedge b) \in Q$, for every $Q \in \tau$. Thus, $a, b \in T$ iff $(a \wedge b) \in T$. \\
To see item (b), let $a \in T$ or $b \in T$, i.e., $a \in Q$, for every $Q \in \tau$ or $b \in Q$, for every $Q \in \tau$. So we have that $a \in Q$ or $b \in Q$, for every $Q \in \tau$. Consequently, by definition \ref{def1}, $(a \vee b) \in Q$, for every $Q \in \tau$, this is, $(a \vee b) \in T$.  Clearly, if $T$ is prime, then the other implication is valid.\\
Let us prove item (c). Let be $a \in T$ and $(a \to b) \in T$, i.e., $a \in Q$, for every $Q \in \tau$ and $(a \to b) \in Q$, for every $Q \in \tau$. So, we have for any  $Q \in \tau$,
$(a \to b) \in Q$, i.e, $\forall P \supseteq Q$, totally prime, if $a \in P$, then $b \in P$. Because $Q$ is totally prime, and $a \in Q$, we must have that $b \in Q$. Thus, $b \in T$.\\
For item (d), observe that the implication from the right to the left is obvious, for  $TPTh_\Le \subseteq PTh_\Le$. It remains to show the other implication. Let $P \in PTh_\Le$ a prime theory generated by $\xi \subseteq TPTh_\Le$, i.e., $P = \bigcap \xi$. Let $(a \to b) \in P$. So, $(a \to b) \in Q$, for every $Q \in \xi$. By definition of implication in \ref{def1}, we have that for every totally prime theory $S \supseteq Q$, if $a \in S$ then $b \in S$. Let now $R \supseteq Q$ a prime theory, such that $a \in R$.  As $R$ is an intersection of totally prime theories, $a$ pertences to every totally prime theory generating the theory $R$. By definition \ref{def1}, we must have $b \in R$, showing (d). \cuad.

\begin{Df} \label{stable}
Let  $\Le,\Le'$ distributive abstract logics.\\
(a)  A {\bf logic application} is a function  $h:Expr_{\Le}\to Expr_{\Le'}$, satisfying $\{h^{-1}(T')|\ T'\in Th_{\Le'}\}\subseteq Th_{\Le}$. We write simply $h:\Le\to \Le'$. \\
(b) A logic application is {\bf stable} iff $\{h^{-1}(P')|\ P'\in PTh_{\Le'}\} \subseteq PTh_{\Le}$. \\
(c) A logic application is {\bf strongly stable} iff $h$ is stable and for every $P' \in PTh_{\Le'}$,  $P \in PTh_\Le$ such that $h^{-1}(P') \subseteq P$, exists $Q' \in PTh_{\Le'}$ such that $P' \subseteq Q'$ and $h^{-1}(Q') = P$.\\
(d) A logic application is {\bf normal} iff $\{h^{-1}(T')|\ T'\in Th_{\Le'}\}=Th_{\Le}$.
\end{Df}

\section{Duality of abstract logics, a general method}
%

\hp In this section, we will introduce the category $\mathcal{LD}$ of distributive abstract logics, cf. \cite{BL11} and develop a  general method to {\sl switch} from abstract logics to lattices. This easy method will allow us to extend some known duality results for lattices to abstract logics. Even almost in an easy manner we can go from abstract logics to lattices, it is not immediately clear that the duality results, will also hold, because the categories always carry with them morphisms.

Remembering theorem \ref{MinGer}, we will in this article always work with abstract logics which are closed by union of chains. Therefore, we always have a set of generators.



In \cite{BL11}, the authors show the analogue  of Stone-Birkhoff's theorem for abstract distributive logics. We will repeat some stuff in this direction.

\begin{Df}
Let $\Le$ an abstract logic with disjunction and $A \subseteq Expr_\Le$. We say that $A$ is {\bf closed under disjunction} iff for every $a, b \in Expr_\Le$, $a \in A$ and $b \in A$, we always have $(a \vee b) \in A$.
\end{Df}

The proof of the next theorem is an application of Zorn's lemma.

\begin{Teo} \label{StoBir}
Let $\mathcal{L}$ be a distributive abstract logic. Let $T \in Th_{\Le}$ and $S \subseteq Expr_\mathcal{L}$ a non empty set closed by  disjunction satisfying $T \cap S  = \emptyset$.
Then, there exists a prime theory $P \in PTh_{\Le}$ such that $T \subseteq P$ and $P \cap S = \emptyset$. \cuad
\end{Teo}

 \begin {Cor} \label{CorSto}
 Let $\mathcal{L}$ a distributive abstract logic. Let $T \in Th_{\Le}$ and  $a \in Expr_\Le$ such that $a \not \in T$. Then there exists a prime theory $P \in PTh_\Le$ such that $P \supseteq T$ and $a \not \in P$.
 \end{Cor}

\dem : Consider for the proof the following set $\overline{ \{ a \} } := \{ b \in Expr_\Le | \hem b \Vdash a \}$. Observe that $\overline{ \{ a \} }$ is closed by disjunction. In fact, being $c, b \in \overline{ \{ a \} }$,  we have that $c \Vdash a$ and $b \Vdash a$. So, by definition of $\Vdash$ we must have $a \in \bigcap \{ T \in Th_\Le | \hem c \in T \}$ and  $a \in \bigcap \{ T \in Th_\Le | \hem b \in T \}$. Consequently,  $a \in \bigcap \{ T \in Th_\Le | \hem c \in T$ or  $b \in T \}$.  By lemma \ref{Lemprelim}, we have that $a \in  \bigcap \{ T \in Th_\Le | \hem (c \vee  b) \in T \}$, i.e., $(c \vee b) \Vdash a$, showing that $\overline{ \{ a \} }$ is closed by disjunction. Furthermore, we have that $T \cap \overline{ \{ a \}} = \emptyset$. By theorem \ref{StoBir}, we obtain the desired. \cuad

The next result is also easy to prove.

\begin{remark}
Let $\Le$ be a distributive abstract logic. Then $PTh_\Le$ is a generator set for $Th_\Le$.
\end{remark}

 \dem : Let $T \in Th_\Le$. Consider $T' := \bigcap \{ P | \hem P \in PTh_\Le$ and $P \supseteq T \}$. We prove that $T = T'$. Clearly $T\subseteq T'$. Suppose that $T \not = T'$. So, there exists $a \in T' \setminus T$. By Corollary \ref{CorSto}, exists a prime theory $P \in PTh_\Le$ such that $P \supseteq T$ and $a \not \in P$. This is a contradiction, because $T'$ was defined as intersection of all prime theories extending  $T$. \cuad

\begin{Lem}
$\mathcal{LD}$ is in fact a category. \cuad
\end{Lem}

In the following we will introduce the almost trivial method to switch from abstract logics to lattices. For this let $\Le=(Expr_{\Le},Th_{\Le},\mathcal{C}_{\Le})$ be a distributive abstract logic. We introduce in  $\Le$ the following order. Let $ a,b\in Expr_{\Le}$,
\begin{equation} \label{order}
 a\leq b\ \Leftrightarrow\ S_{a}\subseteq S_{b},\hem
{\rm with} \hem S_{a}=\{P\in PTh_{\Le};\ a\in P\}
\end{equation}

It is easy to show that $\leq$ is a partial order. For  antisymmetry we use the last theorem  \ref{StoBir}. Now we have a structure of a distributive lattice.

\begin{Lem} \label{loglat}
$A=(Expr_{\Le},\leq)$ is a distributive (bounded) lattice.
\end{Lem}

\dem : First, we show that we have in fact a lattice. For this we establish that \\
\nhp \hfl $inf\{a,b\}=a\wedge b \,$ and  $\,  sup\{a,b\}=a\vee b.  $ \hfl

Let $P\in S_{a\wedge b}$ then $a\wedge b\in P$, by \ref{def1}, $a\in P$ and $b\in P$. Therefore $S_{a\wedge b}\subseteq S_{a}$ and $S_{a\wedge b}\subseteq S_{b}$. Thus
$a\wedge b\leq a$ and  $a\wedge b\leq b$.

Let $c\in Expr_{\Le}$ such that $c\leq a$ and $c\leq b$.
Let $P\in PTh_{\Le}$ such that $c\in P$. By  \ref{order}, we have that $a\in P$ and $b\in P$. By \ref{def1}, $a\wedge b\in P$. Thus $S_{c}\subseteq S_{a\wedge b}$ and so
$c\leq a\wedge b.$

 $sup\{a,b\}=a\vee b$ is showed analogously, with the only exception that the primeness of the theories will play a crucial rule. So, $A$ is a lattice.

For distributivity it suffices to show that  for $a,b,c\in Expr_{\Le}$, we have that $(a\vee b)\wedge(a\vee c)\leq a\vee(b\wedge c)$.

Let $P\in PTh_{\Le}$ such that $(a\vee b)\wedge(a\vee c)\in P$. Then $(a\vee b)\in P$ and $(a\vee c)\in P$. From $a\vee b\in P$, we have $a\in P$ or $b\in P$. Because of $a\vee c\in P$,  $a\in P$ or $c\in P$.

If  $a\in P$, then $P\in S_{a}\subseteq S_{a\vee (b\wedge c)}$.

If $a\not\in P$, we have that $b\in P$ e $c\in P$, thus $b\wedge c\in P$. So, $P\in S_{b\wedge c}\subseteq S_{a\vee(b\wedge c)}$. Therefore
$(a\vee b)\wedge(a\vee c)\leq a\vee(b\wedge c).$

Remark that $S_{\bot}=\emptyset$ and $\emptyset\subseteq S_{a}$ for every $a\in Expr_{\Le}$, so $\bot\leq a$ for every $a$. From $S_{\top}=PTh_{\Le}$, we have that $S_{a}\subseteq S_{\top}$ for every $a\in Expr_{\Le}$, and $A$ is bounded.\cuad

\begin{remark} \label{latlog}
Let $\Omega=(A,\vee,\wedge,\bot,\top)$ be a distributive bounded lattice, then we construct the following distributive abstract logic  $\Le=(Expr_{\Le},Th_{\Le},\mathcal{C}_{\Le})$, with $Expr_{\Le}=A$, $Th_{\Le}=\{F | \hem F$ is a proper filter of $A\}$ and $\mathcal{C}_{\Le}=\{\vee,\wedge,\bot,\top\}$ and $TPTh_\Le := \{ Q | \hem Q$ is a completely prime filter in $A \}$.
\end{remark}

\begin{Lem} \label{HeyLog}
With the above notations, $\Le$ is a distributive abstract logic.
\end{Lem}

\dem : Let $\mathcal{T}\subseteq Th_{\Le}$, then clearly $\bigcap \mathcal{T}$ is a proper filter of $A$. Therefore, $\bigcap\mathcal{T} \in Th_\Le$. It is also clear that this logic is closed by union of chains. The properties for the connectives  follow easily from the filter properties. The distributivity of the logic follows easily from the conditions of  \ref{def1}. Knowing that $\bot\not\in T\ \forall\ T\in Th_{\Le}$ and $\top\in T\ \forall\ T\in Th_{\Le}$, we finish this proof.
 \cuad

\begin{remark}
Let $\Le$ be the distributive abstract logic introduced above. Then $CTh_{\Le}= PTh_\Le =\{ P | \hem P$ is a prime filter of $A\}$.
\end{remark}

\dem : Let $T\in CTh_{\Le}$, so $T\in Th_{\Le}$ and is a proper filter of $A$. As $T\in CTh_{\Le}$, we have by definition that $a\vee b\in T\Leftrightarrow a\in T$ and $ b\in T$, thus $T$ is prime and  $CTh_{\Le}\subseteq PTh_\Le$.
Consider now $P\in PTh_\Le$, this is $P$ is a filter, and so
\nhp \hfl $a\in P$ and $ b\in P \Leftrightarrow a\wedge b\in P.$\hfl

The fact that $P$ is prime, implies that
\nhp \hfl $a\in P$ or $ b\in P\Leftrightarrow a\vee b\in P.$\hfl

So, $PTh_\Le \subseteq CTh_{\Le}$, finishing this proof. \cuad

We define now the category of distributive abstract logics $\mathcal{LD}$\index{Categoria!$\mathcal{LD}$}  as the category with  objects, being  distributive abstract logics and with morphisms  stable logic maps introduced in \ref{stable}.

\begin{Lem}\label{L1}
Let $=_{\leq}$ be the equality defined by the ordering relation introduced in \ref{order} and   $=_{\Le}$  be the equality meaning logical equivalence in the distributive abstract logic $\Le$, i.e., for all $a, b \in Expr_\Le$, $a =_\Le b$ \hem iff \hem $a \Vdash b$ and $b \Vdash a$.
Then these two equalities coincide.
\end{Lem}

\dem :  $(\Rightarrow)$ Let $a=_{\leq}b$, for some $a, b \in Expr_\Le$. Then clearly, $a\leq b$ and $b\leq a$. By  $a\leq b$, we have that $S_{a}\subseteq S_{b}$.  Clearly, $b\in \bigcap\{P\in PTh_{\Le}|\ a\in P\}$. By  \ref{CorSto}, $PTh_{\Le}$ is a generator set for $Th_{\Le}$, and for one $T\in Th_{\Le}$ with $a\in T$, $T$ is intersection of a subset $\mathcal{G}\subseteq PTh_{\Le}$ containing $a$. Therefore $b\in \bigcap\mathcal{G}$, and we can infer that $b\in T$. Thus $a\Vdash_{\Le}b$.
The other case, $b \leq a$ implies $b\Vdash_{\Le}a$, is treated in the same way.

$(\Leftarrow)$ Suppose that $a\Vdash_{\Le}b$, then $b \in \bigcap\{T\in Th_{\Le}|\ a\in T\}$, and so $b\in T$, for every $T\in Th_{\Le}$ with $a\in T$. Particularly, for every $T\in PTh_{\Le}$ such that $a\in T$. So, $S_{a}\subseteq S_{b}$ and we infer that $a\leq b.$
The other case is treated similarly. \cuad

\begin{Cor}
The relations $\leq$ and $\Vdash_{\Le}$ are the same. \cuad
\end{Cor}

In the next Lemma, we will show that stable logic maps are in fact morphisms of the underlying lattices.

\begin{Lem} \label{3.10}
Logic maps in distributive abstract logics are morphisms of lattices.
\end{Lem}

\dem :  Let $\mathcal{L}$ and $\mathcal{L}'$ abstract logics, and  $h: \mathcal{L} \to \mathcal{L}'$ a stable logic map, cf. \ref{stable} (b). We will show that $h$ is a lattice morphism.

For this let $P'\in PTh_{\Le'}$ such that $h(a\vee b)\in P'$. Thus, $a\vee b\in h^{-1}(P')$. Because  $h$ is stable, we have that of $h^{-1}(P')\in PTh_{\Le}$. Therefore, $a\in h^{-1}(P')$ or $b\in h^{-1}(P')$, and so, $h(a)\in P'$ or $h(b)\in P'$.  Consequently, $h(a)\vee h(b)\in P'$, and therefore,
$ S_{h(a\vee b)}\subseteq S_{h(a)\vee h(b)}$. So $h(a\vee b)\leq h(a)\vee h(b)$. The other inequality is proved in the same manner, so that we have  $h(a\vee b)=_{\leq}h(a)\vee h(b)$. And by the above lemma, $h(a\vee b)=_{\Le'}h(a)\vee h(b)$. \\
Completely analogously we can prove, $h(a\wedge b)=_{\Le'} h(a)\wedge h(b)$.\cuad

The next Lemma establishes that in fact lattice morphisms and stable logic maps are the same.

\begin{Lem}
Lattice morphisms in distributive lattices are  stable logic maps in distributive abstract logics.
\end{Lem}

\dem : Let $A, A'$ be distributive lattices and  $f:A\to A'$ a lattice morphism. By construction, we have that $Th_{\Le}:=\{T | \hem T$ is a filter in $A\}.$ Being $f$ a lattice morphism, $f^{-1}(T)$ is a filter in $A$. \\
 In the case of $P'\in PTh_{\Le'}:=\{P' | \hem P'$ is a prime filter in $A'\}$, we have also that, $f^{-1}(P')$ is a prime filter. Thus, $f$ is in fact a stable logic map. \cuad

With the above result, we see that distributive lattices and abstract distributive logics are in bijective correspondence, established in the following way:

(a) Let $\mathcal{L}$ be an abstract distributive logic. Then, we construct a distributive lattice in the way described in \ref{loglat}, denoted by $^{*}\Le := (Expr_\Le, \leq)$, which is in fact a distributive bounded lattice with $inf(a;b) = a \wedge b$ and $sup(a;b)= a \vee b$. \\
This done, we apply the construction in \ref{latlog} to $^{*}\Le$, and we denote by $ _{*}{^{*}\Le}$ the abstract distributive logic obtained. It is not difficult to prove that $ _{*}{^{*}\Le}=\Le$.

Remark that for a theory $T\in Th_{\Le}$,  $T$ is a filter in $^{*}\Le$. Let $a,b\in T$. Because $PTh_{\Le}$ is a generator set for the logic $\Le$, $T=\bigcap\mathcal{P}$, for some $\mathcal{P}\subseteq PTh_{\Le}$.  So $a,b\in \bigcap\mathcal{P}$, and therefore, $a,b\in P$, for every $P\in \mathcal{P}$. \\
But  $P$ is a prime theory, and so $a\wedge b\in P\ \forall\ P\in \mathcal{P}$ and $a\wedge b\in\bigcap\mathcal{P}=T$. \\
Considering $a\in T$ and $b\in Expr_{\Le}$ such that $a\leq b$. By the Lemma \ref{L1}, this is $a\Vdash b$, and so, $b\in \bigcap\{Q\in Th_{\Le}|\ a\in Q\}$. It follows that $b\in T$. Thus, $T$ is a filter of $^{*}\Le.$ \\
Because the theories in $ _{*}{^{*}\Le}$ are proper filters in $^{*}\Le$, we have the desired.

(b) On the other side, we construct by the method \ref{latlog} from a distributive bounded lattice  $A$,  an abstract distributive logic $_{*}A$ and by \ref{loglat} the distributive bounded lattice $^{*}\ _{*} A$. Again, it is not difficult to prove that $^{*}\ _{*} A = A$.

A simple exercise shows that for $\Le,\Le'$ abstract distributive logics and $h:\Le\to\Le'$ a stable logic map, we have that $_{*}\ ^{*}h=h$. Analogously for  $f:A\to A'$ a lattice  morphism, we have that $^{*}\ _{*}f=f$.

So, using the natural transformations, the identity maps, we have established the following theorem.

\begin{Teo}
The category $\mathcal{LD}$ is equivalent with the category $\mathcal{D}ist$\index{Categoria!$\mathcal{D}ist$}, where $\mathcal{D}ist$ denotes the category of distributive lattices. \cuad
\end{Teo}

\begin{Cor}
The category $\mathcal{LD}$ is dually equivalent with the category $\mathcal{P}riest$\index{Categoria!$\mathcal{P}riest$}, where $\mathcal{P}riest$ denotes the category of Priestley spaces. \cuad
\end{Cor}

\begin{Cor}
The category $\mathcal{LD}$ is dually equivalent with the category $\mathcal{S}pec$\index{Categoria!$\mathcal{S}pec$}, where $\mathcal{S}pec$ denotes the category of spectral spaces. \cuad
\end{Cor}

\begin{Cor}
All bitopological dualities, as noted in \cite{bezgabkur} are valid for the category $\mathcal{LD}$. \cuad
\end{Cor}

\section{The case of abstract intuitionistic logics}

\hp In this section, we consider the abstract logic $\mathcal{L}:=(Expr_{\Le},Th_{\Le},\mathcal{C})$, with $\mathcal{C}=\{\wedge,\vee,\rightarrow, \bot, \top \}$, i.e., an  intuitionistic abstract logic in the sense of definition \ref{def1}. We have already constructed a distributive lattice, using the abstract connectives $\mathcal{C}=\{\wedge,\vee, \bot, \top \}$. In the following we will extend these ideas in our new set of connectives. So, we begin to introduce a new category named $\mathcal{LI}$, the category of intuitionistic abstract logics. Then we will show that an intuitionistic abstract logic is an Heyting algebra, and vice versa.

We introduce first the category of intuitionstic abstract logics, whose objects are intuitionistic abstract logics and the morphisms are strongly stable logic maps, defined in definition \ref{stable}. We denote by $\mathcal{LI}$ the category of intuitionistic abstract logic.

\begin{Lem}
$\mathcal{LI}$ is in fact a category.
\end{Lem}

\dem :  Clearly the identiy map is strongly stable and so a morphism in the category. It remains to show that these morphisms are closed under composition. Let $\Le, \Le', \Le'' \in ob(\mathcal{LI})$ and $h: \Le \to \Le'$ e $g:\Le' \to \Le''$ morphisms. We will prove that $g \circ h : \Le \to \Le''$ is a morphism.

Because $h$ is stable, we have that $h^{-1}(P') \in PTh_\Le$, for every $P' \in PTh_{\Le'}$. Using the fact that $g$ is stable,  $h^{-1}(g^{-1}(P'')) \in PTh_\Le$, for every $P'' \in PTh_{\Le''}$, and so $g \circ h$ is also stable.

It remains to show the second condition of definition \ref{stable} (c). For this, let $P'' \in PTh_{\Le''}$ and $P \in PTh_\Le$ be such that $h^{-1}(g^{-1}(P'')) \subseteq P$. Because $g$ is a   morphism, $g^{-1}(P'') \in PTh_{\Le'}$, this is, there exists $P' \in PTh_{\Le'}$ such that $g^{-1}(P'') = P'$. Because  $h$ is a morphism, there is $Q' \in PTh_{\Le'}$ such that $P' \subseteq Q'$ and $P = h^{-1} (Q')$. Thus, $g^{-1}(P'') \subseteq  Q'$. Using the property that  $g$ is strongly stable, there is $Q'' \in PTh_{\Le''}$ such that $P'' \subseteq Q''$ and $Q' = g^{-1}(Q'')$.

So, we have that $P'' \subseteq Q''$ and $P = h^{-1}(Q') = h^{-1}(g^{-1}(Q''))$, finishing the proof that $g\circ  h$ is strongly stable. In fact, $\mathcal{LI}$ forms a category. \cuad

We want to show that $(Expr_\Le; \leq)$ is a Heyting algebra. First the following Lemma.

\begin{Lem}\label{3.1}
Our implication in $\mathcal{C}$ satisfies adjunction, that is, given $z,a,b\in Expr_{\Le}$,
\[z\leq a\rightarrow b \hem \Leftrightarrow \hem z\wedge a\leq b.\]
\end{Lem}

\dem : $(\Rightarrow)$ Suppose that  $z\leq a\rightarrow b$, then $S_{z}\subseteq S_{a\rightarrow b}$.\\
Take $P\in PTh_{\Le}$ such that $z\wedge a\in P$, so $z\in P$ and $a\in P$. Because $z\in P$, we have that $a\rightarrow b\in P$. From the Lemma  \ref{Lemprelim} it follows that
\[a\wedge(a\rightarrow b)\in P\Rightarrow b\in P.\]
Thus, $S_{z\wedge a}\subseteq S_{b}\Rightarrow z\wedge a\leq b.$
\vtres

$(\Leftarrow)$ Let now $z\wedge a\leq b$. Take $P\in PTh_{\Le}$ such that $z\in P$. Let $P'\in PTh_{\Le}$ such that $P\subseteq P'$ and $a\in P'$.\\
Because $z\in P\subseteq P'$, we have that $z\in P'$ and $ a\in P'$, this is, $z\wedge a\in P'$. Thus, $b\in P'$. From $z\wedge a\leq b$, we follow that $a\rightarrow b\in P$. So, $S_{z}\subseteq S_{a\rightarrow b}$, and thus,  $z\leq a\rightarrow b$.\cuad

\begin{Cor} \label{loghey}
$(Expr_\Le; \leq)$ is a Heyting algebra. \cuad
\end{Cor}

Next, we want to show that every Heyting algebra originates  an intuitionistic abstract logic. From the preceding section, we know that every distributive lattice is also a distributive abstract logic, and so the following lemma is sufficient for establishing an intuitionistic abstract logic from every Heyting algebra.

\begin{Lem} \label{heylog}
Let $A$ be a  Heyting algebra. Define $\mathcal{L}:=(Expr_\Le; Th_\Le; \mathcal{C})$ exactly as in remark \ref{latlog}, with the only exception that $\mathcal{C}:= \{ \wedge, \vee, \to, \bot, \top \}$.\\
Then $\mathcal{L}$ is an intuitionistic abstract logic.
\end{Lem}

\dem : It suffices to show that the implication $\to$ satisfies
 the following modified condition of definition \ref{def1}:  for every $a,b\in A$ and prime filter $T$ of $A$,
\[a\rightarrow b\in T\Leftrightarrow for\ every\ prime \ filter \ T'\supseteq T,\ a\in T'\Rightarrow b\in T'\]

 $(\Rightarrow)$ Suppose that $a\rightarrow b\in T$ with $T$ prime filter of $A$. Take $T'$ prime filter such that $T\subseteq T'$ and $a\in T'$. Thus $a\in T'$ and $a\rightarrow b\in T'$, and so $a\wedge (a\rightarrow b)\in T'$. By $a\wedge (a\rightarrow b)\leq b$, we have that $b\in T'$.

$(\Leftarrow)$ Let $T$ be a prime filter and suppose that $a\rightarrow b\not\in T$. Observe that $T \cup \{ a \}$ has the fip (finite intersection property). In the other case, there would be $t_1, \ldots, t_n \in T$ such that $t_1 \wedge \ldots \wedge t_n \wedge a = \bot$. By adjunction, we would have $t_1 \wedge \ldots \wedge t_n \leq a \to \bot$. Because  $T$ is a filter, we have that $t_1\wedge \ldots \wedge t_n \in T$ and  consequently, $(a \to \bot) \in T$. By $\bot \leq b$ and the fact that implication $\to$ is a monotone map, we infer that $a\to \bot \leq a\to b$, and so, $(a \to b) \in T$, a contradiction. Take now  $T\cup\{a\}$ and consider the filter generated $\langle T\cup\{a\}\rangle$, which is proper. We extend this filter to a prime  filter $T'$. Observe that $b\not\in T'$, because in the other case, $b\in T'$, we would have $z\in T$ such that $z\wedge a\leq b$, and by adjunction, $z\leq a\rightarrow b$ and so, once again $a\rightarrow b\in T$, a contradiction, finishing our proof.\cuad

\begin{remark}
Let $\Le$ be an intuitionistic abstract logic. Then $CTh_\Le = PTh_\Le$. \cuad
\end{remark}

\begin{Lem}
Let $h: \Le \to \Le'$ be an intuitionistic logic map, cf.  \ref{stable} (c). Then, $h$ is a morphism of Heyting algebras.
\end{Lem}

\dem : By lemma \ref{3.10}, the morphism $h$ preserves $\wedge$ and $\vee$ and so preserves order. By $S_a \cap S_{a \to b} \subseteq S_b$, we know that $a \wedge (a \to b) \le b$, for $a, b \in Expr_\Le$. Therefore, we have\\
\nhp \hfl $h(a \wedge (a \to b)) = h(a) \wedge' h(a\to b) \leq h(b)$. \hfl

 \nhp By adjunction, cf. lemma \ref{3.1}, we infer that  $h(a \to b) \leq h(a) \to' h(b)$.

 It remains to show that\\
\nhp \hfl $h(a) \to' h(b) \leq h(a \to b)$, this is, $S_{h(a) \to' h(b)} \subseteq S_{h(a \to b)}$. \hfl

\nhp  Let $P' \in PTh_{\Le'}$ such that $P' \not \in S_{h(a \to b)}$, i.e., $h(a \to b) \not \in P'$.

Thus $(a \to b) \not \in h^{-1}(P')$. Observe that $h^{-1}(P') \in PTh_\Le$. By the definition of implication in abstract logics, there exists $P \in PTh_{\Le}$ such that $h^{-1}(P') \subseteq P$ with $a \in P$ and $b \not \in P$. From the second property of definition \ref{stable} (c), there exists  $Q' \in PTh_{\Le'}$ with $P' \subseteq Q'$ and $P = h^{-1}(Q')$.
Thus, $h(a) \in Q'$ and $h(b) \not \in Q'$, i.e., $(h(a) \to' h(b)) \not \in P'$, finishing the proof. \cuad

By the results so far,  given an intuitionistic abstract logic $\mathcal{L}$ we obtain by $^{*}\Le := (Expr_\Le, \le)$ the Heyting algebra using \ref{3.1} and \ref{loghey}. This done we apply the construction \ref{heylog} to $^{*}\Le$ and we denote by $ _{*}{^{*}\Le}$ the intuitionistic abstract logic obtained.\\
 On the other hand using \ref{heylog} and \ref{loghey} we obtain for any Heyting algebra $A$, that $^{*}\ _{*} A = A$.


\begin{Lem} \label{4.7}
Let $h: A \to A'$ a Heyting algebra morphism. With the above notations, $h: \Le \to \Le'$ is an intuitionistic logic map.
\end{Lem}

\dem : It is easy to show that $h$ is a stable logic map. It remains to show the second property of the definition \ref{stable} (c). Let $P \in PTh_\Le := \{ P | \hem P$ is a prime filter in $A \}$ and $P' \in PTh_{\Le'}$ be such that $h^{-1}(P') \subseteq P$. We have to exhibit a prime filter $Q' \in PTh_{\Le'}$ such that $P' \subseteq Q'$ and $P = h^{-1}(Q')$.  Remembering the definition of Esakia morphism and the proof, that every Heyting algebra morphism induces an Esakia morphism, we apply the same proof and obtain the affirmation of our proposition. For the interested reader we give a sketch of this proof in the following remark. \cuad

\begin{remark}
We give a sketch of proofs to be made for finishing the last proposition, see also \cite{dar} and\cite{mor}. \\
(i) Let $A$ be a Heyting algebra, and for $B \subseteq A$, let $\downarrow B = \{ x \in A | \hem \exists y \in B, x \leq y \}$. Denoting $X:= \{ P  | \hem P$ prime filter in $A \}$, we show that $\downarrow (S_a \cap X \setminus S_b) = X \setminus S_{a \to b}$, for all $a,b \in A$. \\
(ii) This done we show that if $Y \subseteq X$ is a clopen subset of $X$ with respect to the Esakia topology, then $Y$ has the form $S_a \cap X \setminus S_b$, for some $a, b \in A$. This fact comes from the compactness of the Esakia space $X$. \\
(iii) In the third step, take a clopen subset $V$ of $X$ such that $P \in V$. By (ii), there exist $a, b \in A$ such that $V = S_a \cap X \setminus S_b$. By $S_a \to S_b = X \setminus \downarrow (S_a \setminus S_b)$ we can show that $h^{-1}(\downarrow V) = \downarrow h^{-1}(V)$. \\
(iv) Now introduce $X:= \{ P | \hem P$ is prime filter in $A \}$ and $X':= \{ P' | \hem P$ is prime filter in $A' \}$. Topologize the two spaces by the Esakia topology and define $h_* : X' \to X$, by $h_*(P') := h^{-1}(P')$. Then we are able to show that $h_*$ is in fact an Esakia morphism, and so particularly, we have that there exists $Q' \in X'$ such that $P' \subseteq Q'$ and $h^{-1}(Q') = P$, finishing the proof of Lemma \ref{4.7}. \cuad

\end{remark}

Denoting the categories of Heyting algebras with the respective Heyting algebra morphisms and of Esakia spaces with the respective Esakia morphisms, cf. \cite{dar} for example, by $\mathcal{H}ey$ and $\mathcal{E}sa$, respectively, we have proved  the following theorem.

\begin{Teo}
The categories $\mathcal{LI}$ and $\mathcal{H}ey$ are dually equivalent. \cuad
\end{Teo}

An immediate corollary, using the known Esakia duality is

\begin{Cor}
The categories $\mathcal{LI}$ and $\mathcal{E}sa$ are dually equivalent. \cuad
\end{Cor}

\end{document}